# LOOP STATISTICS IN THE TOROIDAL HONEYCOMB DIMER MODEL

By Cédric Boutillier and Béatrice de Tilière[1]

*UPMC University Paris 06 and CNRS and Université de Neuchâtel*

The dimer model on a graph embedded in the torus can be interpreted as a collection of random self-avoiding loops. In this paper, we consider the uniform toroidal honeycomb dimer model. We prove that when the mesh of the graph tends to zero and the aspect of the torus is fixed, the winding number of the collection of loops converges in law to a two-dimensional discrete Gaussian distribution. This is known to physicists in more generality from their analysis of toroidal two-dimensional critical loop models and their mapping to the massless free field on the torus. This paper contains the first mathematical proof of this more general physics result in the specific case of the loop model induced by a toroidal dimer model.

**1. Introduction.** Two-dimensional critical loop models are believed by physicists to renormalize at criticality onto a Gaussian free field theory (Coulomb gas). On this basis, they are able to derive explicit formulae for partition functions on the torus [1, 3, 6, 7, 15], from which they obtain information about the asymptotic distribution of the sum of the winding numbers of the loops. In this paper, we give the first mathematical proof of this result, in the case of the toroidal uniform dimer model on the honeycomb lattice.

Let $G = (V(G), E(G))$ be a graph embedded in the torus. A *dimer configuration* or *perfect matching* $M$ of $G$ is a subset of edges of $G$ such that every vertex of $G$ is incident to exactly one edge of $M$. Superimposing $M$ onto a reference dimer configuration $M_0$ yields a collection of self-avoiding loops together with doubled edges; see Figure 2. Loops in this collection can wind around the torus horizontally and vertically.

We consider the case where $G$ is a quotient of the hexagonal lattice and where dimer configurations are chosen uniformly among all dimer configurations of $G$. The algebraic sum of the winding numbers of the loops in

Received September 2006; revised January 2008.
[1]Supported by Swiss National Fund Grant 4 710 2009
*AMS 2000 subject classifications.* 60K35, 82B20.
*Key words and phrases.* Dimers, winding number, loop ensemble.

the collection is then a random variable. The main result of this paper can loosely be stated as follows (see Theorem 2 of Section 1.3 for a precise statement):

THEOREM 1. *When the mesh of the graph $G$ tends to $0$ and the aspect of the torus is fixed, the probability that the algebraic sum of the winding numbers of the loops has horizontal component $k$ and vertical component $\ell$ converges to*

$$\frac{1}{Z_\rho} e^{-\pi(k^2/\rho + \rho\ell^2)/2},$$

*where $Z_\rho = \sum_{(k,\ell)\in\mathbb{Z}^2} e^{-\pi(k^2/\rho+\rho\ell^2)/2}$ is the normalizing factor and $\rho$ is the ratio of the two side lengths of the torus.*

Let us now describe the setting in more detail.

1.1. *Toroidal honeycomb dimer model.* The *dimer model* is a statistical mechanics model introduced to represent the adsorption of diatomic molecules on the surface of a crystal. When, in addition, the underlying graph is bipartite (as is the case for the square lattice or the honeycomb lattice), this model can be interpreted as a random interface model in dimension $2+1$, via the *height function* [17]. The dimer model has the attractive feature of being exactly solvable [8, 9, 16] and is believed to be conformally invariant in the scaling limit; for rigorous results on this, see [5, 10, 11, 12].

In the case where the surface of the crystal is modeled by the regular hexagonal lattice $H$, we speak of the *honeycomb dimer model*. The honeycomb lattice $H$ has a natural embedding in the plane, in which all faces are regular hexagons of side length 1 (or, equivalently, the dual faces are equilateral triangles of side length $\sqrt{3}$). Consider the two vectors $x$ and $y$ represented in Figure 1. The lattice $H$ and its bipartite coloring are invariant under the action of $x$ and $y$ by translation. If, for every $m, n \in \mathbb{N}^*$, we define $\mathbb{L}_{m,n}$ to be the group of translations spanned by $mx$ and $ny$, then $\mathbb{L}_{m,n}$ is a subgroup of the symmetry group of $H$. The toroidal graph $H_{m,n}$ is defined to be the quotient $H_{m,n} = H/\mathbb{L}_{m,n}$.

The graph $H_{m,n}$ can be obtained by cutting out an $mx \times my$ rectangle in $H$ and then gluing opposite sides together. A horizontal side of the rectangle, oriented from left to right, is mapped to an oriented closed loop on the torus, as is a vertical side of the rectangle, oriented from bottom to top. Let us denote by $\gamma^h$ and $\gamma^v$ these two oriented closed loops. The first homology group $H_1(\mathbb{T})$ is spanned by the homology classes of $\gamma^h$ and $\gamma^v$.

The *modulus* of the torus is the ratio of the two complex numbers representing the vertical and the horizontal sides of the rectangle. Since the length of the horizontal side is $3m$ and the vertical side has length $n\sqrt{3}$ by construction, the modulus of $H_{m,n}$ is $i\frac{n}{\sqrt{3}m}$.

A *dimer configuration* of $H_{m,n}$ is a perfect matching of $H_{m,n}$, that is, a subset of edges $M$ of $H_{m,n}$ such that every vertex of $H_{m,n}$ is incident to exactly one edge of $M$. Let us denote by $\mathcal{M}(H_{m,n})$ the set of dimer configurations of $H_{m,n}$. Suppose that a positive weight function $\nu$ is assigned to edges of $H_{m,n}$, that is, every edge $e$ has weight $\nu(e)$. Then, every dimer configuration $M$ of $H_{m,n}$ has an energy $\mathcal{E}(M) = -\sum_{e \in M} \log \nu(e)$. The probability of occurrence of the dimer configuration $M$ of $H_{m,n}$ is given by the Boltzmann measure $\mu_{m,n}$:

$$\mu_{m,n}(M) = \frac{e^{-\mathcal{E}(M)}}{Z_{m,n}(\nu)} = \frac{\prod_{e \in M} \nu(e)}{Z_{m,n}(\nu)},$$

where $Z_{m,n}(\nu) = \sum_{M \in \mathcal{M}(H_{m,n})} \prod_{e \in M} \nu(e)$ is the normalizing constant, known as the *partition function*.

When the weight function $\nu$ is periodic, it is known, using classical subadditivity arguments (as in [4]) that the quantity $-\frac{1}{mn} \log Z_{m,n}(\nu)$ converges when $m$ and $n$ tend to $\infty$. The limit is denoted by $\mathsf{f}(\nu)$ and is called the *free energy per fundamental domain*.

In this paper, we consider dimer configurations of $H_{m,n}$ chosen with respect to the uniform measure (i.e., the Boltzmann measure corresponding to weights 1 on all edges). The corresponding partition function is called the *uniform partition function*. In this case, the free energy per fundamental domain, simply denoted by $\mathsf{f}$, is given by the following formula [8], whose derivation is recalled in Section 3.2:

$$(1) \qquad \mathsf{f} = -\frac{1}{4\pi^2} \int_0^{2\pi} \int_0^{2\pi} \log(2(\cos\psi + 1) - e^{i\phi}) \, d\phi \, d\psi.$$

1.2. *Toroidal dimer model and self-avoiding loops.* The dimer model on the toroidal graph $H_{m,n}$ can be interpreted as a collection of self-avoiding loops as follows. Let $M_0$ be a fixed dimer configuration of $H_{m,n}$ and let $M$ be any other dimer configuration of $H_{m,n}$. The superimposition of $M_0$ and $M$ then consists of self-avoiding *doubled edges* and *alternating loops*, where *doubled edges* are edges covered by a dimer in both $M_0$ and $M$, and *alternating loops* are cycles whose edges are, in alternation, dimers of $M_0$ and $M$; see Figure 2. This feature is due to the fact that, by definition of perfect matchings, every vertex of $H_{m,n}$ is incident to exactly one edge of $M_0$ and one edge of $M$.

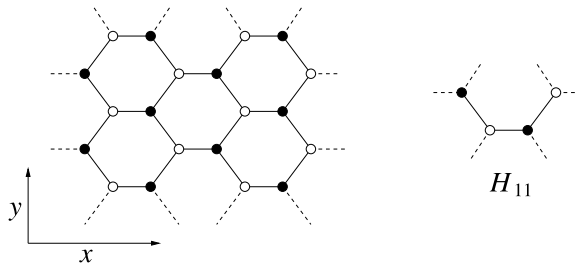

FIG. 1. *The lattice $H$ is invariant under the action of $x$ and $y$ by translation (left). Fundamental domain $H_{1,1}$ (right).*

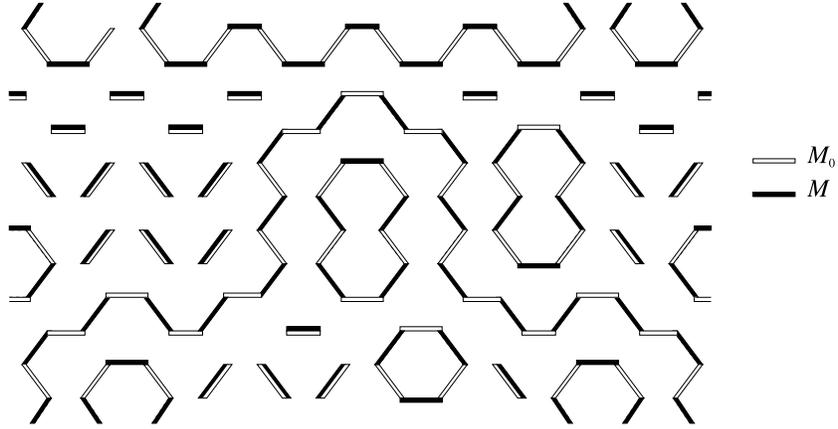

FIG. 2. *The superimposition of $M$ and $M_0$ consists of doubled edges and alternating loops.*

Orienting the dimers of $M$ from their white end to their black end and the dimers of $M_0$ from their black end to their white end gives rise to an orientation of the loops. Let us denote by $M \ominus M_0$ the set of oriented loops obtained from this superimposition; see Figure 3 for an example.

A loop $C$ of $M \ominus M_0$ can then be seen as a closed path on the torus $\mathbb{T}$. The equivalence class $[C]$ of $C$ in the first homology group $H_1(\mathbb{T}) \simeq \mathbb{Z}^2$ can be decomposed in the basis $([\gamma^h], [\gamma^v])$. Its coordinates $(C^h, C^v)$ in this basis are called the *winding number of $C$*: $C^h$ (resp., $C^v$) is the algebraic number of times the loop $C$ winds horizontally (resp., vertically) around the torus. The *winding number of the dimer configuration $M$*, denoted by $\text{wind}_{M_0}(M)$, is the sum of the winding numbers of all of the loops contained in $M \ominus M_0$:

$$\text{wind}_{M_0}(M) = \sum_{\substack{C \text{ loop} \\ \text{in } M \ominus M_0}} [C] \in \mathbb{Z}^2.$$

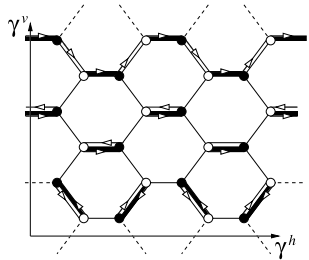

FIG. 3. *The superimposition $M \ominus M_0$ consists of oriented loops. In this example, the winding number is* $\mathrm{wind}_{M_0}(M) = (1,0)$.

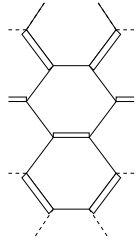

Fig. 4. *Dimer configuration of $H_{1,3}$ which generates the reference matching $M_0$.*

An example of a computation of $\text{wind}_{M_0}(M)$ is given in Figure 3. Note that the dependence on $M_0$ of the winding number of $M$ is quite simple: if $M_1$ is another dimer configuration, then

$$\text{wind}_{M_1}(M) = \text{wind}_{M_0}(M) - \text{wind}_{M_0}(M_1).$$

From now on, we assume that $n$ is a multiple of 3 and fix the reference dimer configuration $M_0$ of $H_{m,n}$ to be the one generated by translations of the dimer configuration of $H_{1,3}$ of Figure 4. We also drop the subscript $M_0$ in $\text{wind}_{M_0}(M)$.

Since dimer configurations of $H_{m,n}$ are chosen according to the uniform measure, for every $m, n$, $\text{wind}(\cdot)$ is a random variable. Let us call it the *winding number* and, in order to stress the dependence on $m$ and $n$, denote it by $\text{wind}_{m,n}(\cdot)$.

The main result of this paper is an explicit expression for the asymptotic distribution of the random variables $(\text{wind}_{m,n})$ when $m, n$ tend to infinity and when the modulus $i\frac{n}{\sqrt{3}m}$ of $H_{m,n}$ converges to $i\rho$ for some $\rho > 0$.

1.3. *Statement of result.* Let us assume that the modulus $i\frac{n}{\sqrt{3}m}$ of $H_{m,n}$ converges to $i\rho$ for some $\rho > 0$. Recall that when $n$ is a multiple of 3, $\text{wind}_{m,n}$ is the winding number of the uniformly distributed dimer configurations of $H_{m,n}$, computed with respect to the reference dimer configuration $M_0$ defined above. Then, the main result of this paper is the following theorem:

THEOREM 2. *In the joint limit $m, n \to \infty$, $\frac{n}{\sqrt{3}m} \to \rho$, the sequence of random variables $(\mathrm{wind}_{m,n})$ converges in distribution to the two-dimensional discrete Gaussian random variable $\mathrm{wind}_\rho$ whose law is given by*

$$(2) \qquad \forall (k, \ell) \in \mathbb{Z}^2 \qquad \mathbb{P}[\mathrm{wind}_\rho = (k, \ell)] = \frac{1}{Z_\rho} e^{-\pi(k^2/\rho + \rho \ell^2)/2},$$

*where $Z_\rho = \sum_{(k,\ell) \in \mathbb{Z}^2} e^{-\pi(k^2/\rho + \rho \ell^2)/2}$.*

- A similar result was obtained by Kenyon and Wilson [14] in the case of the square lattice embedded in the *cylinder*. Working on the torus makes computations much more difficult since it means dealing with the toroidal partition function in the proof, which is a combination of four terms [9] (instead of one, as in the cylinder case). Moreover, we have to extract information about the two components of the winding number (instead of one, as in the cylinder case). Note, also, that in proving Theorem 2, we give a *full asymptotic expansion* of a perturbation of the uniform partition function; see Theorem 4, Section 2 below.
- It was brought to our attention, after the acceptance of this paper, that an asymptotic expansion of the uniform partition function was obtained in the physics literature by Ferdinand, in the case of the square lattice [7]. Nevertheless, let us stress the following facts: the expansion of [7] is not perturbative, it is not done to the same level of mathematical rigor and no information about the distribution of the winding number is inferred.
- Theorem 2 holds if the honeycomb lattice is replaced by the square lattice. We believe that the techniques and ideas applied here could be used to extend the result to quotients of the honeycomb and square lattice, embedded on tori with modulus in the upper half complex plane (and not just $i\rho$, $\rho > 0$). More generally, we conjecture the result to be true when considering the dimer model on any periodic bipartite graph within the liquid phase [13], with an appropriate embedding.

1.4. *Outline of the paper.*

- In Section 2, we prove that the moment generating function of $\mathrm{wind}_{m,n}$ can be expressed in terms of a perturbed uniform partition function $Z_{m,n}(\alpha, \beta)$; see Section 2 for a definition. Theorem 4 then gives a full asymptotic expansion of $Z_{m,n}(\alpha, \beta)$, from which we deduce pointwise convergence of the moment generating function of $\mathrm{wind}_{m,n}$ to the moment generating function of $\mathrm{wind}_\rho$ and, hence, Theorem 2.
- The remainder of the paper consists of the proof of Theorem 4, giving the full asymptotic expansion of the perturbed uniform partition function $Z_{m,n}(\alpha, \beta)$. For the reader's convenience, the proof is split into two parts:

- By [8], the partition function $Z_{m,n}(\alpha,\beta)$ can be expressed as a linear combination of four terms, $Z_{m,n}^{(\sigma\eta)}(\alpha,\beta)$, $\sigma, \eta \in \{0,1\}$. Proposition 8 of Section 3 gives a full asymptotic expansion for each of the four terms $Z_{m,n}^{(\sigma\eta)}(\alpha,\beta)$ as a function of Jacobi's four elliptic theta functions. The proof of Proposition 8 is postponed until Section 4. Proposition 9 then gives an explicit expression for the combination of Jacobi theta functions involved in the expression of $Z_{m,n}(\alpha,\beta)$. The proof of Theorem 4 is thus completed, apart from the proof of Proposition 8.
- Section 4 consists of the proof of Proposition 8, giving the full asymptotic expansions of the four terms $Z_{m,n}^{(\sigma\eta)}(\alpha,\beta)$ as functions of Jacobi's four elliptic theta functions.

**2. Winding number and partition function.** Define $F_{m,n}$ to be the moment generating function of the random variable $\text{wind}_{m,n}$:

$$\forall (\alpha,\beta) \in \mathbb{R}^2 \qquad F_{m,n}(\alpha,\beta) = \mathbb{E}[e^{-\pi \text{wind}_{m,n} \cdot (\alpha,\beta)}]$$
$$= \sum_{(k,\ell) \in \mathbb{Z}^2} \mathbb{P}[\text{wind}_{m,n} = (k,\ell)] e^{-\pi(\alpha k + \beta \ell)}.$$

Convergence in distribution of the sequence $(\text{wind}_{m,n})$ to the two-dimensional Gaussian random variable $\text{wind}_\rho$ given in (2) is equivalent to pointwise convergence of the sequence $(F_{m,n})$ to the corresponding moment generating function.

Lemma 3 below gives an expression of $F_{m,n}$ in terms of a perturbed uniform partition function, defined as follows. Let us introduce the appropriate choice of perturbed edge-weights. Define edges to be of type I (resp., II, III), as in Figure 5. Then, for $\alpha, \beta \in \mathbb{R}$, let us assign weights $a = e^{-\alpha\pi/(2m)}$ to edges of type I, $b^{-1} = e^{-\beta\pi/(2n)}$ to edges of type II and $b = e^{\beta\pi/(2n)}$ to edges of type III. In the sequel, we shall also use the notation $A = a^{-2m} = e^{\alpha\pi}$, $B = b^{2n} = e^{\beta\pi}$. Observe that the weights $a, b^{-1}, b$ tend to 1 when $m, n$ tend to infinity and thus yield a perturbation of the uniform partition function. They are used to collect information on the uniform measure. Let us denote by $Z_{m,n}(\alpha,\beta)$ the partition function of the graph $H_{m,n}$ corresponding to these weights and let us call it the *perturbed uniform partition function* or, in short, *perturbed partition function*.

LEMMA 3. *The moment generating function $F_{m,n}(\alpha,\beta)$ of $\text{wind}_{m,n}$ and the perturbed uniform partition function $Z_{m,n}(\alpha,\beta)$ are related in the following way:*

$$F_{m,n}(\alpha,\beta) = e^{\pi\alpha n/3} \frac{Z_{m,n}(\alpha,\beta)}{Z_{m,n}(0,0)}.$$

PROOF. For $i =$I, II, III, let $N_i(M)$ be the number of edges of type $i$ in the dimer configuration $M$. Then, by definition, the partition function $Z_{m,n}(\alpha,\beta)$ is

$$Z_{m,n}(\alpha,\beta) = \sum_{M\in\mathcal{M}(H_{m,n})} (e^{-\alpha\pi/(2m)})^{N_{\mathrm{I}}(M)}(e^{-\beta\pi/(2n)})^{N_{\mathrm{II}}(M)}(e^{\beta\pi/(2n)})^{N_{\mathrm{III}}(M)}$$

$$= \sum_{M\in\mathcal{M}(H_{m,n})} A^{-N_{\mathrm{I}}(M)/(2m)} B^{(N_{\mathrm{III}}(M)-N_{\mathrm{II}}(M))/(2n)}$$

where, we recall, $A = e^{\pi\alpha}, B = e^{\pi\beta}$.

Let us compute the two components $\mathrm{wind}^h_{m,n}$ and $\mathrm{wind}^v_{m,n}$ of $\mathrm{wind}_{m,n}$ as functions of $N_{\mathrm{I}}, N_{\mathrm{II}}, N_{\mathrm{III}}$. Consider the $2n$ left-to-right horizontal paths of the dual graph $H^*_{m,n}$, as in Figure 6 (left). Then, for any such path $\gamma$, $\mathrm{wind}^v_{m,n}(M)$ is equal to the number of positive (right-to-left) crossings of $M \ominus M_0$ along $\gamma$ minus the number of negative (left-to-right) crossings of $M \ominus M_0$ along $\gamma$. Summing over all $2n$ horizontal paths and observing that $N_{\mathrm{III}}(M_0) = N_{\mathrm{II}}(M_0)$, we obtain

$$2n \cdot \mathrm{wind}^v_{m,n}(M) = N_{\mathrm{II}}(M) + N_{\mathrm{III}}(M_0) - N_{\mathrm{III}}(M) - N_{\mathrm{II}}(M_0)$$
$$= N_{\mathrm{II}}(M) - N_{\mathrm{III}}(M).$$

In a similar way, considering the $2m$ top-to-bottom vertical paths of Figure 6 (right) and observing that $N_{\mathrm{I}}(M_0) = 2mn/3$, we obtain

$$2m \cdot \mathrm{wind}^h_{m,n}(M) = N_{\mathrm{I}}(M) - N_{\mathrm{I}}(M_0) = N_{\mathrm{I}}(M) - \frac{2mn}{3}.$$

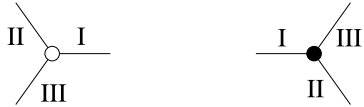

FIG. 5. *The three types of edges around a white vertex (left) and around a black vertex (right).*

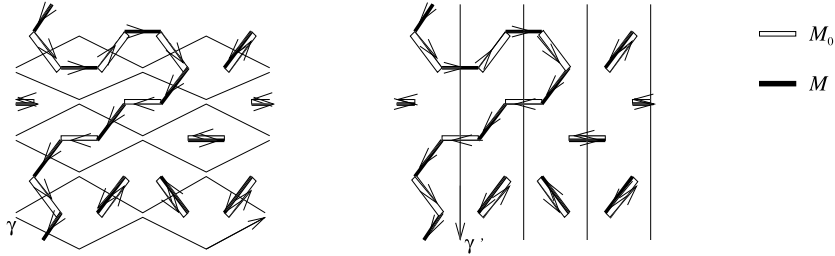

FIG. 6. *Computing $\mathrm{wind}^v_{m,n}$ (left), $\mathrm{wind}^h_{m,n}$ (right).*

Plugging this into the partition function $Z_{m,n}(\alpha,\beta)$, we obtain

$$Z_{m,n}(\alpha,\beta) = A^{-n/3} \sum_{M\in\mathcal{M}(H_{m,n})} A^{-\text{wind}^h_{m,n}(M)} B^{-\text{wind}^v_{m,n}(M)}$$

$$= e^{-\pi\alpha n/3} \sum_{(k,\ell)\in\mathbb{Z}^2} C_{k,\ell} e^{-\pi(\alpha k + \beta \ell)},$$

where $C_{k,\ell}$ is the number of dimer configurations whose winding number is $(k,\ell)$. The proof is completed by recalling that dimer configurations of $H_{m,n}$ are chosen with respect to the uniform measure, which implies that

$$\mathbb{P}[\text{wind}_{m,n} = (k,\ell)] = \frac{C_{k,\ell}}{Z_{m,n}(0,0)}$$

$[Z_{m,n}(0,0)$ is the uniform partition function$]$.

□

Theorem 4 below gives a precise asymptotic expansion of $Z_{m,n}(\alpha,\beta)$. Combined with Lemma 3, this yields pointwise convergence of the sequence $(F_{m,n})$ to the moment generating function $F_\rho$ of the discrete Gaussian random variable $\text{wind}_\rho$ and, hence, Theorem 2.

THEOREM 4. *In the joint limit $m,n \to \infty$, $\frac{n}{\sqrt{3}m} \to \rho$, we have the following asymptotic expansion for the perturbed partition function $Z_{m,n}(\alpha,\beta)$. For all $(\alpha,\beta) \in \mathbb{R}^2$,*

$$Z_{m,n}(\alpha,\beta) = (-1)^{mn} e^{-\pi n\alpha/3} e^{-mn\mathsf{f}} \frac{e^{\pi\rho/6}}{\sqrt{2\rho}P(e^{-\rho\pi})^2}$$

$$\times \sum_{(k,\ell)\in\mathbb{Z}^2} e^{\pi(\alpha k + \beta \ell)} e^{-\pi(k^2/\rho + \rho\ell^2)/2}(1+o(1)),$$

*where $P(q) = \prod_{k=1}^\infty (1-q^{2k})$ and $\mathsf{f}$ is the free energy per fundamental domain of equation (1).*

REMARK 5. Recall that, by definition, the free energy per fundamental domain $\mathsf{f}$ is given by

$$\mathsf{f} = -\lim_{m,n\to\infty} \frac{1}{mn} \log Z_{m,n}(0,0).$$

Hence, it is not surprising that $\mathsf{f}$ should govern the exponential growth rate $mn$ of $Z_{m,n}(\alpha,\beta)$.

COROLLARY 6. *For all $(\alpha, \beta) \in \mathbb{R}^2$,*
$$\lim_{\substack{m,n \to \infty \\ n/\sqrt{3}m \to \rho}} F_{m,n}(\alpha, \beta) = F_\rho(\alpha, \beta).$$

PROOF. Combining Lemma 3 and Theorem 4, we have
$$\lim_{\substack{m,n \to \infty \\ n/(\sqrt{3}m) \to \rho}} F_{m,n}(\alpha, \beta) = \frac{1}{Z_\rho} \sum_{(k,\ell) \in \mathbb{Z}^2} e^{\pi(\alpha k + \beta \ell)} e^{-\pi(k^2/\rho + \rho \ell^2)/2}$$
$$= \frac{1}{Z_\rho} \sum_{(k,\ell) \in \mathbb{Z}^2} e^{-\pi(\alpha k + \beta \ell)} e^{-\pi(k^2/\rho + \rho \ell^2)/2}$$
$$= F_\rho(\alpha, \beta),$$
where $Z_\rho = \sum_{(k,\ell) \in \mathbb{Z}^2} e^{-\pi(k^2/\rho + \rho \ell^2)/2}$. In the second equality, we have used the fact that the sum is symmetric in $(k, \ell)$. □

REMARK 7. Recall that $F_{m,n}(\alpha, \beta)$ is the moment generating function of the random variable $\text{wind}_{m,n}$, which is computed using the reference dimer configuration $M_0$. A natural question which arises is what happens when the reference dimer configuration is changed. Let $F_{m,n}^{M_1}(\alpha, \beta)$ be the moment generating function of the winding number, computed using a generic reference matching $M_1$. Then, looking at the proof of Lemma 3, we obtain
$$F_{m,n}^{M_1}(\alpha, \beta) = A^{N_{\text{I}}(M_1)/(2m)} B^{(-N_{\text{III}}(M_1) + N_{\text{II}}(M_1))/(2n)} \frac{Z_{m,n}(\alpha, \beta)}{Z_{m,n}(0,0)}.$$
Using Theorem 4, we deduce that
$$F_{m,n}^{M_1}(\alpha, \beta) = A^{N_{\text{I}}(M_1)/(2m) - n/3} B^{(-N_{\text{III}}(M_1) + N_{\text{II}}(M_1))/(2n)}$$
$$\times \frac{1}{Z_\rho} \sum_{(k,\ell) \in \mathbb{Z}^2} e^{\pi(\alpha k + \beta \ell)} e^{-\pi(k^2/\rho + \rho \ell^2)/2} (1 + o(1)).$$
Hence, Corollary 6 holds provided that the reference dimer configuration satisfies $N_{\text{I}}(M_1) = \frac{2mn}{3}$ and $N_{\text{III}}(M_1) = N_{\text{II}}(M_1)$. Observing that $N_{\text{I}}, N_{\text{II}}, N_{\text{III}}$ are always constrained to satisfy $N_{\text{I}} + N_{\text{II}} + N_{\text{III}} = 2mn$, this implies that $N_{\text{I}}(M_1) = N_{\text{II}}(M_1) = N_{\text{III}}(M_1) = \frac{2mn}{3}$.

**3. Proof of Theorem 4.** The main ingredient in the proof of Theorem 4 is an explicit expression for the perturbed partition function $Z_{m,n}(\alpha, \beta)$ as a combination of Jacobi theta functions. This is given in Proposition 8, Section 3.2 below. The proof of Proposition 8 is postponed until Section 4. In Section 3.1, we recall the definition of Jacobi theta functions. Proposition 9 of Section 3.3 gives a concise formula for the combination of theta functions involved in the expression of the perturbed partition function $Z_{m,n}(\alpha, \beta)$. Section 3.4 consists of the proof of Theorem 4, using all of the above.

3.1. *Jacobi theta functions.* Recall the definition of Jacobi's four elliptic theta functions $\vartheta_i(\zeta, q)$, $i = 1, \ldots, 4$, and their expressions in terms of infinite products:

$$\vartheta_1(\zeta, q) = \sum_{k=-\infty}^{\infty} (-1)^{k-1/2} e^{(2k+1)i\zeta} q^{(k+1/2)^2}$$

$$= 2q^{1/4} \sin(\zeta) P(q) \prod_{\ell=1}^{\infty} (1 - 2q^{2\ell} \cos(2\zeta) + q^{4\ell}),$$

$$\vartheta_2(\zeta, q) = \sum_{k=-\infty}^{\infty} e^{(2k+1)i\zeta} q^{(k+1/2)^2}$$

$$= 2q^{1/4} \cos(\zeta) P(q) \prod_{\ell=1}^{\infty} (1 + 2q^{2\ell} \cos(2\zeta) + q^{4\ell}),$$

$$\vartheta_3(\zeta, q) = \sum_{k=-\infty}^{\infty} e^{2ki\zeta} q^{k^2} = P(q) \prod_{\ell=0}^{\infty} (1 + 2q^{2\ell+1} \cos(2\zeta) + q^{4\ell+2}),$$

$$\vartheta_4(\zeta, q) = \sum_{k=-\infty}^{\infty} (-1)^k e^{2ki\zeta} q^{k^2} = P(q) \prod_{\ell=0}^{\infty} (1 - 2q^{2\ell+1} \cos(2\zeta) + q^{4\ell+2}),$$

where $P(q) = \prod_{k=1}^{\infty}(1 - q^{2k})$. It is sometimes convenient to use the notation $\vartheta_i(\zeta|\tau)$ for $\vartheta_i(\zeta, e^{i\pi\tau})$.

The proofs of these formulas, as well as other properties of Jacobi theta functions, can be found in [2].

3.2. *Perturbed uniform partition function.* Let us recall the definition of the perturbed partition function $Z_{m,n}(\alpha, \beta)$. It is the partition function of the dimer model on the graph $H_{m,n}$, where edges are assigned perturbed uniform weights: edges of type I (resp., II, III) have weights $a = e^{-\alpha\pi/(2m)}$ (resp., $b^{-1} = e^{-\beta\pi/(2n)}$, $b$). Recall, also, the notation $A = e^{\alpha\pi}$, $B = e^{\beta\pi}$ and $\rho = \lim_{m,n\to\infty} \frac{n}{\sqrt{3}m}$. An explicit formula for the partition function $Z_{m,n}(\alpha, \beta)$ is given in [8]:

(3) $$Z_{m,n}(\alpha, \beta) = \tfrac{1}{2}((-1)^n(-Z_{m,n}^{(00)}(\alpha, \beta) + Z_{m,n}^{(01)}(\alpha, \beta)) \\ + Z_{m,n}^{(10)}(\alpha, \beta) + Z_{m,n}^{(11)}(\alpha, \beta)),$$

where

$$Z_{m,n}^{(\sigma\eta)}(\alpha, \beta) = \prod_{z^m=(-1)^\sigma} \prod_{w^n=(-1)^\eta} P(z, w),$$

$$(4) \quad P(z,w) = \det \begin{pmatrix} \dfrac{1}{b} + \dfrac{b}{w} & a \\ az & b + \dfrac{w}{b} \end{pmatrix}$$

$$= \left(\frac{1}{b} + \frac{b}{w}\right)\left(b + \frac{w}{b}\right) - a^2 z = \frac{w}{b^2} + \frac{b^2}{w} + 2 - a^2 z.$$

Each of the four terms $Z_{m,n}^{(\sigma\eta)}(\alpha,\beta)$ is a determinant of a *Kasteleyn matrix*, a relative of the adjacency matrix of the graph $H_{m,n}$, computed with discrete Fourier transforms, using the invariance of the graph under translations by the vectors $x$ and $y$ of Figure 1. Using Riemann sums and analyzing the behavior of $P(z,w)$ in the neighborhood of its zeros, one deduces (see, e.g., Kasteleyn [8]) that the free energy per fundamental domain with weights 1 on the edges is

$$f = -\lim_{m,n\to\infty} \frac{1}{mn} \log Z_{m,n}(0,0) = \frac{1}{4\pi^2} \int_{\mathbb{T}} \log\left(w + \frac{1}{w} + 2 - z\right) \frac{dz}{z} \frac{dw}{w},$$

$$= -\frac{1}{4\pi^2} \int_0^{2\pi} \int_0^{2\pi} \log(2(\cos\psi + 1) - e^{i\phi}) \, d\phi \, d\psi,$$

which is precisely equation (1).

The following proposition gives the asymptotic expansion of the four terms $Z_{m,n}^{(\sigma\eta)}(\alpha,\beta)$ involved in the explicit expression (3) of $Z_{m,n}(\alpha,\beta)$. The proof is postponed until Section 4.

PROPOSITION 8. *In the joint limit $m,n \to \infty$, $\frac{n}{\sqrt{3}m} \to \rho$, we have the following asymptotic expansion for the four terms involved in the partition function $Z_{m,n}(\alpha,\beta)$:*

$$\left.\begin{array}{r}(-1)^n Z_{m,n}^{(00)}(\alpha,\beta) \\ (-1)^n Z_{m,n}^{(01)}(\alpha,\beta) \\ Z_{m,n}^{(11)}(\alpha,\beta) \\ Z_{m,n}^{(10)}(\alpha,\beta)\end{array}\right\} = (-1)^{mn} A^{-n/3} e^{-mnf} \frac{e^{\pi\alpha^2\rho/2} e^{\pi\rho/6}}{P(q)^2}$$

$$\times \begin{cases} (-\vartheta_1(\zeta,q)\vartheta_1(\bar{\zeta},q)) + o(1), \\ \vartheta_2(\zeta,q)\vartheta_2(\bar{\zeta},q) + o(1), \\ \vartheta_3(\zeta,q)\vartheta_3(\bar{\zeta},q) + o(1), \\ \vartheta_4(\zeta,q)\vartheta_4(\bar{\zeta},q) + o(1), \end{cases}$$

*where $\zeta = \frac{\pi}{2}(\rho\alpha + i\beta)$, $q = e^{-\rho\pi}$ and $f$ is given in equation (1).*

3.3. *Recombining Jacobi theta functions.* The following proposition gives an explicit expression for the combination of Jacobi theta functions involved in the perturbed partition function.

PROPOSITION 9.

$$\forall \zeta = x + iy \in \mathbb{C}, \forall \tau \in \mathbb{H} = \{z \in \mathbb{C}; \operatorname{Im} z > 0\},$$

$$\sum_{i=1}^{4} \vartheta_i(\zeta|\tau)\vartheta_i(\bar{\zeta}|\tau) = \sqrt{\frac{2i}{\tau}} e^{-i2x^2/(\pi\tau)} \vartheta_3\left(\frac{x}{\tau}\bigg|-\frac{1}{2\tau}\right)\vartheta_3\left(iy\bigg|\frac{\tau}{2}\right).$$

PROOF. Let $q = e^{i\pi\tau}$. By definition of Jacobi theta functions, we have

$$\vartheta_1(\zeta|\tau)\vartheta_1(\bar{\zeta}|\tau) + \vartheta_2(\zeta|\tau)\vartheta_2(\bar{\zeta}|\tau)$$
$$= \sum_{(k,\ell)\in\mathbb{Z}^2} (1-(-1)^{k+\ell}) e^{i(2k+1)\zeta} e^{i(2\ell+1)\bar{\zeta}} q^{(k+1/2)^2+(\ell+1/2)^2}$$
$$= \sum_{(k,\ell)\in\mathbb{Z}^2} (1-(-1)^{k+\ell}) e^{i2(k+\ell+1)x} e^{i2(k-\ell)iy} q^{(k+1/2)^2+(\ell+1/2)^2}.$$

The general term of this sum is nonzero if and only if $k \not\equiv \ell \bmod 2$. In this case, we introduce the two integers $u$ and $v$ such that $k+\ell+1 = 2u$ and $k - \ell = 2v+1$, and rewrite $\vartheta_1(\zeta|\tau)\vartheta_1(\bar{\zeta}|\tau) + \vartheta_2(\zeta|\tau)\vartheta_2(\bar{\zeta}|\tau)$ as

$$\vartheta_1(\zeta|\tau)\vartheta_1(\bar{\zeta}|\tau) + \vartheta_2(\zeta|\tau)\vartheta_2(\bar{\zeta}|\tau)$$
(5)
$$= 2 \sum_{(u,v)\in\mathbb{Z}^2} e^{i2u(2x)} e^{i2(2v+1)iy} q^{2u^2+(2v+1)^2/2}.$$

Similarly, we write

$$\vartheta_3(\zeta|\tau)\vartheta_3(\bar{\zeta}|\tau) + \vartheta_4(\zeta|\tau)\vartheta_4(\bar{\zeta}|\tau)$$
$$= \sum_{(k,\ell)\in\mathbb{Z}^2} (1+(-1)^{k+\ell}) e^{i2k\zeta} e^{i2\ell\bar{\zeta}} q^{k^2+\ell^2}$$
$$= \sum_{(k,\ell)\in\mathbb{Z}^2} (1+(-1)^{k+\ell}) e^{i2(k+\ell)x} e^{i2(k-\ell)iy} q^{k^2+\ell^2}.$$

In this case, the general term is nonzero when $k \equiv \ell \bmod 2$. Setting $u = \frac{k+\ell}{2}$ and $v = \frac{k-\ell}{2}$, we get

(6)
$$\vartheta_3(\zeta|\tau)\vartheta_3(\bar{\zeta}|\tau) + \vartheta_4(\zeta|\tau)\vartheta_4(\bar{\zeta}|\tau)$$
$$= 2 \sum_{(u,v)\in\mathbb{Z}^2} e^{i2u(2x)} e^{i2(2v)iy} q^{2u^2+(2v)^2/2}.$$

Summing (5) and (6), and recombining the sum over $v$, we obtain

$$\sum_{i=1}^{4} \vartheta_i(\zeta|\tau)\vartheta_i(\bar{\zeta}|\tau) = 2 \sum_{(u,v)\in\mathbb{Z}^2} e^{i2u(2x)} e^{i2v(iy)} q^{2u^2+v^2/2}$$

$$= 2\left(\sum_{u=-\infty}^{\infty} e^{2iu(2x)}(q^2)^{u^2}\right)\left(\sum_{v=-\infty}^{\infty} u^{2iv(iy)}(q^{1/2})^{v^2}\right)$$

$$= 2\vartheta_3(2x|2\tau)\vartheta_3\left(iy\left|\frac{\tau}{2}\right.\right)$$

$$= \sqrt{\frac{2i}{\tau}}e^{-i2x^2/(\pi\tau)}\vartheta_3\left(\frac{x}{\tau}\left|-\frac{1}{2\tau}\right.\right)\vartheta_3\left(iy\left|\frac{\tau}{2}\right.\right).$$

The last equality results from the Jacobi identity which describes the transformation of the function $\vartheta_3$ under the modular group [2]:

$$\forall u \in \mathbb{C}, \forall \sigma \in \mathbb{H} \qquad \vartheta_3(u|\sigma) = \sqrt{\frac{i}{\sigma}}e^{-iu^2/(\pi\sigma)}\vartheta_3\left(\frac{u}{\sigma}\left|-\frac{1}{\sigma}\right.\right),$$

with $u = 2x$ and $\sigma = 2\tau$. □

COROLLARY 10. *In the case where* $\zeta = \frac{\pi}{2}(\rho\alpha + i\beta)$ *and* $\tau = i\rho$ *(i.e.,* $q = e^{-\pi\rho}$*),*

$$\sum_{i=1}^{4} \vartheta_i(\zeta,q)\vartheta_i(\bar\zeta,q) = \sqrt{\frac{2}{\rho}}e^{-\pi\alpha^2\rho/2}\vartheta_3\left(-\frac{i\pi\alpha}{2}\left|\frac{i}{2\rho}\right.\right)\vartheta_3\left(\frac{i\pi\beta}{2}\left|\frac{i\rho}{2}\right.\right)$$

$$= \sqrt{\frac{2}{\rho}}e^{-\pi\alpha^2\rho/2}\sum_{(k,\ell)\in\mathbb{Z}^2} A^k B^\ell e^{-\pi(k^2/\rho+\rho\ell^2)/2},$$

*where* $A = e^{\pi\alpha}$ *and* $B = e^{\pi\beta}$.

3.4. *Proof of Theorem 4.* Using expression (3) for the perturbed partition function $Z_{m,n}(\alpha,\beta)$ and Proposition 8, we know that in the joint limit $m, n \to \infty$, $\frac{n}{\sqrt{3}m} \to \rho$,

(7)
$$Z_{m,n}(\alpha,\beta) = (-1)^{mn} A^{-n/3} e^{-mn\mathsf{f}} \frac{e^{\pi\alpha^2\rho/2}e^{\pi\rho/6}}{2P(q)^2}$$
$$\times \left(\sum_{i=1}^{4}\vartheta_i(\zeta,q)\vartheta_i(\bar\zeta,q) + o(1)\right).$$

Replacing the combination of Jacobi theta functions of (7) by the expression given in Corollary 10 yields

(8)
$$Z_{m,n}(\alpha,\beta) = (-1)^{mn} A^{-n/3} e^{-mn\mathsf{f}} \frac{e^{\pi\rho/6}}{\sqrt{2\rho}P(q)^2}$$
$$\times \sum_{(k,\ell)\in\mathbb{Z}^2} A^k B^\ell e^{-\pi(k^2/\rho+\rho\ell^2)/2}(1+o(1)),$$

which is precisely Theorem 4.

**4. Asymptotic expansion of the perturbed partition function: Proof of Proposition 8.** Let us quickly recall the content of Proposition 8. It gives explicit expressions for the four terms $Z_{m,n}^{(00)}(\alpha,\beta)$, $Z_{m,n}^{(01)}(\alpha,\beta)$, $Z_{m,n}^{(11)}(\alpha,\beta)$, $Z_{m,n}^{(10)}(\alpha,\beta)$ involved in the expression (3) of the perturbed partition function $Z_{m,n}(\alpha,\beta)$, as functions of Jacobi's four theta functions. Sections 4.1–4.4 consist of preliminary computations for each of the four terms $Z_{m,n}^{(\sigma\eta)}(\alpha,\beta)$. Proposition 8 is then proved in Section 4.5 using the aforementioned computations. Let us recall the notation used: $a = e^{-\pi\alpha/(2m)}$, $b = e^{\pi\beta/(2n)}$, $A = e^{\pi\alpha}$, $B = e^{\pi\beta}$, $\rho = \lim_{m,n\to\infty} \frac{n}{\sqrt{3}m}$, $\zeta = \frac{\pi}{2}(\rho\alpha + i\beta)$ and $q = e^{-\rho\pi}$.

4.1. *Computations for $Z_{m,n}^{(11)}(\alpha,\beta)$.* The beginning of this computation is inspired by [14]. Using the explicit expression (3) for $Z_{m,n}^{(11)}(\alpha,\beta)$ yields

$$Z_{m,n}^{(11)}(\alpha,\beta) = \prod_{z^m=-1}\prod_{w^n=-1}\left(\frac{w}{b^2} + (2-a^2z) + \frac{b^2}{w}\right)$$

$$= \prod_{z^m=-1}\prod_{w^n=-1}\frac{b^2}{w}\left(\frac{w^2}{b^4} + \frac{w}{b^2}(2-a^2z) + 1\right)$$

$$= \prod_{z^m=-1}\prod_{w^n=-1}\frac{b^2}{w}\left(\frac{w}{b^2} - r_1\right)\left(\frac{w}{b^2} - r_2\right)$$

$$= \prod_{z^m=-1}\prod_{w^n=-1}\frac{1}{b^2w}(w - b^2r_1)(w - b^2r_2),$$

where

$$r_1, r_2 = -1 + \frac{a^2z}{2} \pm \sqrt{\left(1-\frac{a^2z}{2}\right)^2 - 1}.$$

Since $r_1$ and $r_2$ depend on $a$ and $z$ through the product $a^2z$, we define $\phi \in \mathbb{C}$ as $a^2z = e^{i\phi}$ and

(9) $$r_{1,2}(\phi) = -1 + \frac{e^{i\phi}}{2} \pm i\sqrt{1 - \left(1-\frac{e^{i\phi}}{2}\right)^2}.$$

Recall that $a = e^{-\pi\alpha/(2m)}$ and that $z$ is an $m$th root of $-1$. Hence, for large $m$ and $n$, we are concerned with values of $\phi$ whose imaginary part is close to 0 and whose real part is in $[-\pi, \pi]$.

*Choice of $r_1, r_2$ and their domain of definition.* In the formula for $r_1(\phi), r_2(\phi)$, there is an ambiguity in the determination of the square root, which we now clarify. The functions $r_1(\phi), r_2(\phi)$ are analytic on a domain $D$, provided the square root does not vanish in $D$. The square root vanishes when $(\text{Re}(\phi), \text{Im}(\phi)) = (0[2\pi], \log(1/4))$, so let us choose

$$D = \{\phi \in \mathbb{C} | -\pi < \text{Re}(\phi) < \pi, -c < \text{Im}(\phi) < c\},$$

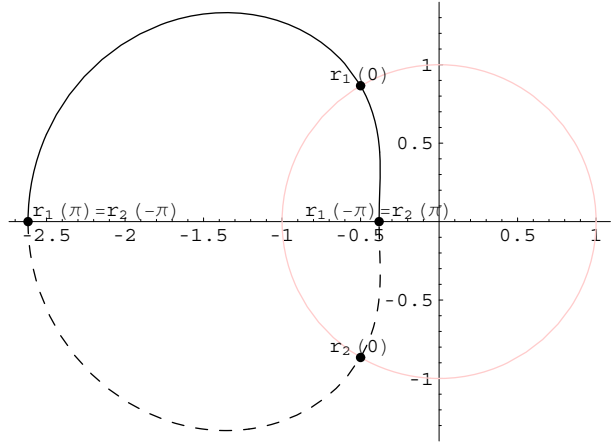

Fig. 7. *Plot of the two functions $r_1(\phi)$ (plain line) and $r_2(\phi)$ (dashed line) for $\phi \in [-\pi, \pi]$, and their position with respect to the unit circle. Note that when $\alpha$ is not 0 and $m$ is large enough, $r_1(i\frac{\pi\alpha}{m})$ and $r_2(i\frac{\pi\alpha}{m})$ are still on the unit circle (see Lemma 11 below).*

where $c$ is some positive constant smaller than $\log(4)$. Now, observe that $r_1(\phi)r_2(\phi) = 1$ so that one of $r_1$, $r_2$ has nonnegative imaginary part. Moreover, $r_1$ satisfies the polynomial equation

$$(10) \qquad r_1 + \frac{1}{r_1} = e^{i\phi} - 2.$$

Taking imaginary parts on both sides yields

$$(11) \qquad \operatorname{Im}(r_1(\phi))\left(1 - \frac{1}{|r_1(\phi)|^2}\right) = e^{-\operatorname{Im}(\phi)} \sin(\operatorname{Re}(\phi)).$$

Hence, $\operatorname{Im}(r_1(\phi)) \neq 0$ whenever $\operatorname{Re}(\phi) \notin \{-\pi, 0, \pi\}$. When $\operatorname{Re}(\phi) = 0$, we look at equation (9) and deduce that the term in the square root is positive when $|\operatorname{Im}(\phi)| < c$ so that $\operatorname{Im}(r_1(\phi)) \neq 0$ in $D$. Using the analyticity of $r_1$ and $r_2$ in $D$, we conclude that one of $r_1$, $r_2$ has imaginary part which is positive in $D$. By convention, we choose $r_1(\phi)$ such that $\operatorname{Im}(r_1(\phi)) > 0$ in $D$. The values for $r_1$, $r_2$ on the boundary of the domain $D$ are chosen by continuity. Figure 7 shows a plot of $r_1(\phi)$, $r_2(\phi)$ when $\phi$ is real.

Let us first compute the product over $w$ in (8). Since $w$ is an $n$th root of $-1$, we have

$$\forall \lambda \in \mathbb{C} \qquad \prod_{w^n=-1}(\lambda - w) = \lambda^n + 1.$$

In particular,

$$\prod_{w^n=-1} w = (-1)^n, \qquad \prod_{w^n=-1}(b^2 r_j(\phi) - w) = 1 + B r_j^n(\phi), \qquad j = 1, 2,$$

so that

$$Z^{(11)}_{m,n}(\alpha,\beta) = (-1)^{mn} B^{-m} \prod_{z^m=-1} (1+Br_1^n(\phi))(1+Br_2^n(\phi)). \tag{12}$$

*Ideas governing the next steps of the computation.* There is a change of behavior in the product (12) when $|r_j(\phi)|$ is smaller, greater or close to 1. We start by factoring out $r_j(\phi)$ when $|r_j(\phi)| > 1$. This part of $Z^{(11)}_{m,n}(\alpha,\beta)$ grows exponentially, so we expect it to involve the free energy per fundamental domain f. We postpone its study until Section 4.5. Then, in Lemma 12, we deal with the part of $Z^{(11)}_{m,n}(\alpha,\beta)$ where $|r_j(\phi)| < 1$. This part makes no contribution to $Z^{(11)}_{m,n}(\alpha,\beta)$. Finally, in Lemma 13, we deal with the part of $Z^{(11)}_{m,n}(\alpha,\beta)$ where $|r_j(\phi)|$ is close to 1. This part involves the Jacobi theta functions giving the discrete Gaussian term of $Z^{(11)}_{m,n}(\alpha,\beta)$. In Sections 4.2, 4.3 and 4.4, we study, in a concise way, the computations for $Z^{(10)}_{m,n}(\alpha,\beta)$, $Z^{(01)}_{m,n}(\alpha,\beta)$ and $Z^{(00)}_{m,n}(\alpha,\beta)$.

In order to proceed with these steps, we need:

1. a characterization of $|r_j(\phi)| > 1$, $|r_j(\phi)| < 1$ and $|r_j(\phi)| = 1$;
2. an expansion of $r_j(\phi)$ when $\phi$ is close to 0 (for Lemmas 12, 13);
3. to know that $|r_1(\phi)|$ is increasing as a function of $\text{Re}(\phi)$ when $-\pi \leq \text{Re}(\phi) < 0$ and $\text{Im}(\phi)$ is small (for Lemma 12).

All of this is provided by the following lemma.

LEMMA 11.  *On the domain $\overline{D}$, the functions $r_1$ and $r_2$ have the following properties:*

1. – *If $0 < \text{Re}(\phi) \leq \pi$, then $|r_1(\phi)| > 1$, $|r_2(\phi)| < 1$.*
   – *If $-\pi \leq \text{Re}(\phi) < 0$, then $|r_1(\phi)| < 1$, $|r_2(\phi)| > 1$.*
   – *If $\text{Re}(\phi) = 0$, then $|r_1(\phi)| = |r_2(\phi)| = 1$.*
2. *Their expansions in the vicinity of $\phi = 0$ are given by*

$$r_1(\phi) = e^{i2\pi/3}\left(1 + \frac{\phi}{\sqrt{3}} + O(\phi^2)\right) = e^{i2\pi/3} e^{\phi/\sqrt{3}+O(\phi^2)}, \tag{13}$$

$$r_2(\phi) = e^{-i2\pi/3}\left(1 - \frac{\phi}{\sqrt{3}} + O(\phi^2)\right) = e^{-i2\pi/3} e^{-\phi/\sqrt{3}+O(\phi^2)}. \tag{14}$$

3. *When $-\pi \leq \text{Re}(\phi) < 0$, $|r_1(\phi)|$ is increasing as a function of $\text{Re}(\phi)$.*

PROOF.  1. Since $r_2(\phi)r_1(\phi) = 1$, it suffices to prove the properties for $r_1(\phi)$. Recall that, by convention, $\text{Im}(r_1(\phi)) > 0$ in $D$. Looking at equation (11), we deduce that when $\sin(\text{Re}(\phi)) > 0$ or, equivalently, $\text{Re}(\phi) \in (0,\pi)$, we have $|r_1(\phi)| > 1$. We further deduce that when $\sin(\operatorname{Re}(\phi)) < 0$ or, equivalently, $\operatorname{Re}(\phi) \in (-\pi, 0)$, we have $|r_1(\phi)| < 1$. The explicit evaluation of $r_1(\pm\pi)$ shows that the inequalities are still strict when $\operatorname{Re}(\phi) = \pm\pi$.

We thank the referee for the idea behind the above argument.

Let us now consider the case where $\operatorname{Re}(\phi) = 0$. Since $\phi \in D$, we have $\operatorname{Im}(\phi) > \log(1/4)$ and

$$1 - \left(1 - \frac{e^{-\operatorname{Im}(\phi)}}{2}\right)^2 > 0$$

so that, in this case, $r_1(\phi)$ is given by

$$r_1(\phi) = -1 + \frac{e^{-\operatorname{Im}(\phi)}}{2} + i\sqrt{1 - \left(1 - \frac{e^{-\operatorname{Im}(\phi)}}{2}\right)^2}$$

with no ambiguity in the square root. An explicit computation yields $|r_1(\phi)| = 1$.

2. The second part of the lemma is obtained via an explicit Taylor expansion of $r_1(\phi), r_2(\phi)$ in a neighborhood 0, which can be performed since these functions are analytic.

3. We need to prove that $\frac{\partial |r_1(\phi)|^2}{\partial \operatorname{Re}(\phi)}$ is nonnegative whenever $-\pi \leq \operatorname{Re}(\phi) < 0$. The derivative of $|r_1(\phi)|^2$ with respect to $\operatorname{Re}(\phi)$ is

$$\frac{\partial |r_1(\phi)|^2}{\partial \operatorname{Re}(\phi)} = 2\operatorname{Re}\left(\frac{\partial r_1(\phi)}{\partial \operatorname{Re}(\phi)} \overline{r_1(\phi)}\right).$$

Observe that $r_1(\phi)$ satisfies the polynomial equation

$$r_1(\phi)^2 + r_1(\phi)(2 - e^{i\phi}) + 1 = 0.$$

From this, we deduce, by differentiating with respect to $\operatorname{Re}(\phi)$, that

$$\frac{\partial r_1(\phi)}{\partial \operatorname{Re}(\phi)} = \frac{ie^{i\phi} r_1(\phi)}{2r_1(\phi) + 2 - e^{i\phi}},$$

$$\frac{\partial |r_1(\phi)|^2}{\partial \operatorname{Re}(\phi)} = 2\operatorname{Re}\left(\frac{ie^{i\phi} |r_1(\phi)|^2}{2r_1(\phi) + 2 - e^{i\phi}}\right).$$

Hence,

$$\operatorname{sign}\left(\frac{\partial |r_1(\phi)|^2}{\partial \operatorname{Re}(\phi)}\right) = -\operatorname{sign}\left(\operatorname{Im}\left(\frac{1}{2e^{-i\phi}(r_1(\phi) + 1) - 1}\right)\right)$$

$$= \operatorname{sign}(\operatorname{Im}(e^{-i\phi}(r_1(\phi) + 1))).$$

Rewriting the quadratic equation satisfied by $r_1(\phi)$ as

$$(r_1(\phi) + 1)^2 = e^{i\phi} r_1(\phi),$$

we see that $e^{-i\phi}(r_1(\phi) + 1)$ equals $\frac{r_1(\phi)}{1+r_1(\phi)}$, whose imaginary part has the same sign as $\text{Im}(r_1(\phi))$, which is nonnegative by definition. $\square$

Since $z$ is an $m$th root of $-1$, we can write $z = e^{i\pi(2j+1)/m}$ for $j \in \{-\lfloor\frac{m}{2}\rfloor, \ldots, \lfloor\frac{m-1}{2}\rfloor\}$ so that $\phi = \frac{\pi}{m}(i\alpha + 2j + 1)$. Using Lemma 11, we know that when $j \geq 0$, $|r_1(\phi)| > 1$, and when $j \leq -1$, $|r_2(\phi)| > 1$. Factoring those terms out in (12) and using an algebraic manipulation to produce a product over positive $j$'s only, we obtain the following expression for $Z_{m,n}^{(11)}(\alpha, \beta)$:

$$Z_{m,n}^{(11)}(\alpha, \beta)$$
$$= (-1)^{mn} B^{-m} \prod_{j=-\lfloor m/2 \rfloor}^{\lfloor (m-1)/2 \rfloor} \left(1 + Br_1^n\left(\frac{\pi(i\alpha + 2j + 1)}{m}\right)\right)$$
$$\times \left(1 + Br_2^n\left(\frac{\pi(i\alpha + 2j + 1)}{m}\right)\right)$$
$$= (-1)^{mn} \left(\prod_{j=0}^{\lfloor (m-1)/2 \rfloor} r_1\left(\frac{\pi(i\alpha + 2j + 1)}{m}\right)\right)^n$$
$$\times \left(\prod_{j=0}^{\lfloor m/2 \rfloor - 1} r_2\left(\frac{\pi(i\alpha - (2j + 1))}{m}\right)\right)^n$$
$$\times \prod_{j=0}^{\lfloor (m-1)/2 \rfloor} \left(1 + B^{-1} r_1^{-n}\left(\frac{\pi(i\alpha + 2j + 1)}{m}\right)\right)$$
$$\times \left(1 + Br_2^n\left(\frac{\pi(i\alpha + 2j + 1)}{m}\right)\right)$$
$$\times \prod_{j=0}^{\lfloor m/2 \rfloor - 1} \left(1 + Br_1^n\left(\frac{\pi(i\alpha - (2j + 1))}{m}\right)\right)$$
$$\times \left(1 + B^{-1} r_2^{-n}\left(\frac{\pi(i\alpha - (2j + 1))}{m}\right)\right).$$

LEMMA 12.    *In the joint limit $m, n \to \infty$, $\frac{n}{\sqrt{3}m} \to \rho$, we have*

$$\prod_{j=\lfloor m^{1/4} \rfloor}^{\lfloor (m-1)/2 \rfloor} \left(1 + B^{-1} r_1^{-n}\left(\frac{\pi(i\alpha + 2j + 1)}{m}\right)\right) = 1 + o(1),$$

$$\prod_{j=\lfloor m^{1/4} \rfloor}^{\lfloor (m-1)/2 \rfloor} \left( 1 + Br_2^n \left( \frac{\pi(i\alpha + 2j + 1)}{m} \right) \right) = 1 + o(1),$$

$$\prod_{j=\lfloor m^{1/4} \rfloor}^{\lfloor m/2 \rfloor - 1} \left( 1 + Br_1^n \left( \frac{\pi(i\alpha - (2j+1))}{m} \right) \right) = 1 + o(1),$$

$$\prod_{j=\lfloor m^{1/4} \rfloor}^{\lfloor m/2 \rfloor - 1} \left( 1 + B^{-1} r_2^{-n} \left( \frac{\pi(i\alpha - (2j+1))}{m} \right) \right) = 1 + o(1).$$

PROOF. We give the proof for the first product on the second line. The proof for the other terms is obtained using the same arguments.

Taking the logarithm of the left-hand side, we get

$$\left| \sum_{j=\lfloor m^{1/4} \rfloor}^{\lfloor m/2 \rfloor - 1} \log \left( 1 + Br_1^n \left( \frac{\pi(i\alpha - (2j+1))}{m} \right) \right) \right|$$

$$\leq \sum_{j=\lfloor m^{1/4} \rfloor}^{\lfloor m/2 \rfloor - 1} \left| \log \left( 1 + Br_1^n \left( \frac{\pi(i\alpha - (2j+1))}{m} \right) \right) \right|$$

$$\leq \sum_{j=\lfloor m^{1/4} \rfloor}^{\lfloor m/2 \rfloor - 1} cB \left| r_1^n \left( \frac{\pi(i\alpha - (2j+1))}{m} \right) \right| \qquad \text{for some constant } c > 0$$

$$\leq mcB \max_{\lfloor m^{1/4} \rfloor \leq j \leq \lfloor m/2 \rfloor - 1} \left| r_1^n \left( \frac{\pi(i\alpha - (2j+1))}{m} \right) \right|$$

$$\leq mcB \left| r_1^n \left( \frac{\pi(i\alpha - 2m^{1/4})}{m} \right) \right| \qquad \text{(Lemma 11, part 3)}$$

$$\leq mcB e^{-2\pi\rho m^{1/4}(1+o(1))} \qquad \text{(Lemma 11, part 2)}$$

$$= o(1). \qquad \square$$

LEMMA 13. *In the joint limit $m, n \to \infty$, $\frac{n}{\sqrt{3}m} \to \rho$, we have*

$$\prod_{j=0}^{\lfloor m^{1/4} \rfloor - 1} \left( 1 + B^{-1} r_1^{-n} \left( \frac{\pi(i\alpha + 2j + 1)}{m} \right) \right) \left( 1 + Br_1^n \left( \frac{\pi(i\alpha - (2j+1))}{m} \right) \right)$$

$$= \frac{\vartheta_3(\bar{\zeta}, q)}{P(q)} + o(1),$$

$$\prod_{j=0}^{\lfloor m^{1/4} \rfloor -1} \left(1 + Br_2^n\left(\frac{\pi(i\alpha + 2j + 1)}{m}\right)\right)\left(1 + B^{-1}r_2^{-n}\left(\frac{\pi(i\alpha - (2j+1))}{m}\right)\right)$$

$$= \frac{\vartheta_3(\zeta, q)}{P(q)} + o(1),$$

where $\zeta = \frac{\pi}{2}(\rho\alpha + i\beta)$, $q = e^{-\rho\pi}$ and $P(q) = \prod_{j=1}^{+\infty}(1 - q^{2k})$.

PROOF. The arguments used are similar to those of Lemma 12. First, note that

$$\vartheta_3(\zeta, q) = P(q)\prod_{j=0}^{\infty}(1 + e^{2i\zeta}q^{2j+1})(1 + e^{-2i\zeta}q^{2j+1})$$

$$= \lim_{m \to \infty} P(q) \prod_{j=0}^{\lfloor m^{1/4} \rfloor -1} (1 + e^{2i\zeta}q^{2j+1})(1 + e^{-2i\zeta}q^{2j+1}).$$

Let us prove the following (the other three cases are handled similarly):

(15)
$$\prod_{j=0}^{\lfloor m^{1/4} \rfloor -1} \left(1 + B^{-1}r_2^{-n}\left(\frac{\pi(i\alpha - (2j+1))}{m}\right)\right)$$
$$= \prod_{j=0}^{\lfloor m^{1/4} \rfloor -1} (1 + e^{2i\zeta}q^{2j+1}) + o(1).$$

A generic term of the right-hand side of (15) can be written as

$$(1 + e^{2i\zeta}q^{2j+1}) = (1 + e^{-\beta\pi}e^{i\pi\alpha\rho}e^{-\pi\rho(2j+1)}).$$

Note that such a term can be equal to 0, whenever $\alpha = \frac{2k+1}{\rho}$ for some integer $k$ and $\beta = -\rho(2j+1)$ for some $j \in \{0, \ldots, \lfloor m^{1/4} \rfloor - 1\}$. So, the product of the right-hand side of (15) might vanish. Observe that when this is the case, there is only one term in the product which is actually 0; let us denote by $\overline{j}$ the corresponding index.

We first consider the part of the product corresponding to nonzero terms on the right-hand side of (15). That is, depending on the values of $\alpha$ and $\beta$, we either consider the whole product or the product over all indices except $\overline{j}$. In order to keep the notation simple, it is implicit that $\overline{j}$ is omitted whenever it occurs. Dividing the left-hand side by the right-hand side in (15) and taking logarithms yields

$$\left|\sum_{j=0}^{\lfloor m^{1/4} \rfloor -1} \log\left(\frac{1 + B^{-1}r_2^{-n}(\pi(i\alpha - (2j+1))/m)}{1 + B^{-1}e^{i\pi\alpha\rho}q^{2j+1}}\right)\right|$$

$$\leq \sum_{j=0}^{\lfloor m^{1/4}\rfloor - 1} \left|\log \frac{1 + B^{-1} r_2^{-n}(\pi(i\alpha - (2j+1))/m)}{1 + B^{-1} e^{i\pi\alpha\rho} q^{2j+1}}\right|.$$

Observe that the two quantities $B^{-1} e^{i\pi\alpha\rho} q^{2j+1}$ and $B^{-1} r_2^{-n}(\frac{\pi(i\alpha - (2j+1))}{m})$ are bounded, and bounded away from $-1$, uniformly in $j = 0, \ldots, \lfloor m^{1/4}\rfloor$, $j \neq \overline{j}$, for $m$ and $n$ large enough. We therefore get the following bounds, where $c_1, c_2$ are positive constants:

$$\left|\sum_{j=0}^{\lfloor m^{1/4}\rfloor - 1} \log\left(\frac{1 + B^{-1} r_2^{-n}(\pi(i\alpha - (2j+1))/m)}{1 + B^{-1} e^{i\pi\alpha\rho} q^{2j+1}}\right)\right|$$

$$\leq \sum_{j=0}^{\lfloor m^{1/4}\rfloor - 1} c_1 \cdot B^{-1} \left|r_2^{-n}\left(\frac{i\alpha - (2j+1)}{m}\right) - e^{i\pi\alpha\rho} q^{2j+1}\right|$$

$$\leq \sum_{j=0}^{\lfloor m^{1/4}\rfloor - 1} c_1 \cdot B^{-1} \left|e^{\pi(i\alpha - (2j+1))n(1+o(1))/(\sqrt{3}m)} \right.$$

$$\left. - e^{\pi(i\alpha - (2j+1))\rho}\right| \qquad \text{by (13)}$$

$$\leq c_2 \sum_{j=0}^{\lfloor m^{1/4}\rfloor} |i\alpha - (2j+1)| e^{-\pi(2j+1)(\rho+o(1))} \left|\rho - \frac{n(1+o(1))}{\sqrt{3}m}\right| = o(1)$$

because the series with general term $|i\alpha - (2j+1)| e^{-\pi(2j+1)\rho}$ is convergent and thus the upper bound goes to zero when $\frac{n}{\sqrt{3}m}$ goes to $\rho$.

Let us now consider the case when $\alpha = \frac{2k+1}{\rho}$ for some integer $k$ and $\beta = -\rho(2\overline{j} + 1)$ for some $\overline{j} \in \{0, \ldots, \lfloor m^{1/4}\rfloor - 1\}$. For all terms involving indices $j \neq \overline{j}$, the above estimate holds. When $j = \overline{j}$, the corresponding term on the right-hand side of (15) is 0; let us show that the corresponding term on the left-hand side is $o(1)$:

$$\left|1 + B^{-1} r_2^{-n}\left(\frac{\pi(i\alpha - (2\overline{j} + 1))}{m}\right)\right|$$

$$= |1 + B^{-1} e^{\pi(i\alpha - (2\overline{j}+1))n(1+o(1))/(\sqrt{3}m)}| \qquad \text{by (13)}$$

$$= |1 + B^{-1} e^{\pi(i\alpha - (2\overline{j}+1))\rho} [e^{\pi(i\alpha - (2\overline{j}+1))(n(1+o(1))/(\sqrt{3}m) - \rho)}]|$$

$$= \left|1 - \left[1 + \pi(i\alpha - (2\overline{j}+1))\left(\frac{n(1+o(1))}{\sqrt{3}m} - \rho\right) + o(1)\right]\right|$$

$$= \left|\pi(i\alpha - (2\overline{j}+1))\left(\frac{n(1+o(1))}{\sqrt{3}m} - \rho\right) + o(1)\right| = o(1).$$

In the two last lines, we have used the assumption $1+B^{-1}e^{\pi(i\alpha-(2\bar{j}+1))\rho}=0$. Combining this with the estimate when $j\neq\bar{j}$ yields the lemma. $\square$

COROLLARY 14. *In the joint limit $m,n\to\infty$, $\frac{n}{\sqrt{3}m}\to\rho$, we have*

$$\frac{(-1)^{mn}Z^{(11)}_{m,n}(\alpha,\beta)}{\Lambda^1_{m,n}(\alpha)\Lambda^2_{m,n}(\alpha)} = \frac{\vartheta_3(\zeta,q)\vartheta_3(\bar{\zeta},q)}{P(q)^2} + o(1),$$

*where*

$$\Lambda^1_{m,n}(\alpha) = \left(\prod_{j=0}^{\lfloor(m-1)/2\rfloor} r_1\left(\frac{\pi(i\alpha+2j+1)}{m}\right)\right)^n,$$

$$\Lambda^2_{m,n}(\alpha) = \left(\prod_{j=0}^{\lfloor m/2\rfloor-1} r_2\left(\frac{\pi(i\alpha-(2j+1))}{m}\right)\right)^n.$$

4.2. *Computations for $Z^{(10)}_{m,n}(\alpha,\beta)$.* The computations for $Z^{(10)}_{m,n}(\alpha,\beta)$ go through in a similar way. Let us just stress where the differences occur. Using the explicit expression of (3) for $Z^{(10)}_{m,n}(\alpha,\beta)$ yields

$$Z^{(10)}_{m,n}(\alpha,\beta) = \prod_{z^m=-1}\prod_{w^n=1}\frac{1}{b^2w}(w-b^2r_1)(w-b^2r_2).$$

This time, $w$ is an $n$th root of $1$ (instead of $-1$), so

$$\forall\lambda\in\mathbb{C}\qquad \prod_{w^n=1}(\lambda-w) = \lambda^n - 1.$$

In particular, $\prod_{w^n=1}w=(-1)^{n-1}$ and, when performing the product over $w$, we obtain

$$Z^{(10)}_{m,n}(\alpha,\beta)$$
$$=(-1)^{m(n-1)}B^{-m}$$
$$\times\prod_{j=-\lfloor m/2\rfloor}^{\lfloor(m-1)/2\rfloor}\left(1-Br_1^n\left(\frac{\pi(i\alpha+2j+1)}{m}\right)\right)\left(1-Br_2^n\left(\frac{\pi(i\alpha+2j+1)}{m}\right)\right)$$
$$=(-1)^{mn}\left(\prod_{j=0}^{\lfloor(m-1)/2\rfloor} r_1\left(\frac{\pi(i\alpha+2j+1)}{m}\right)\right)^n$$
$$\times\left(\prod_{j=0}^{\lfloor m/2\rfloor-1} r_2\left(\frac{\pi(i\alpha-(2j+1))}{m}\right)\right)^n$$
$$\times\prod_{j=0}^{\lfloor(m-1)/2\rfloor}\left(1-B^{-1}r_1^{-n}\left(\frac{\pi(i\alpha+2j+1)}{m}\right)\right)$$

$$\times \left(1 - Br_2^n\left(\frac{\pi(i\alpha + 2j + 1)}{m}\right)\right)$$

$$\times \prod_{j=0}^{\lfloor m/2 \rfloor - 1} \left(1 - Br_1^n\left(\frac{\pi(i\alpha - (2j+1))}{m}\right)\right)$$

$$\times \left(1 - B^{-1}r_2^{-n}\left(\frac{\pi(i\alpha - (2j+1))}{m}\right)\right).$$

The rest of the computation goes through in the same way, except that the $+$'s of $Z_{m,n}^{(11)}$ are replaced by the $-$'s of $Z_{m,n}^{(10)}$. As a consequence, the analog of Lemma 13 involves the fourth Jacobi theta function $\vartheta_4(\zeta, q)$ instead of the third one. We summarize the expression for $Z_{m,n}^{(10)}$ in the following result.

COROLLARY 15. *In the joint limit $m, n \to \infty$, $\frac{n}{\sqrt{3}m} \to \rho$,*

$$\frac{(-1)^{mn} Z_{m,n}^{(10)}(\alpha, \beta)}{\Lambda_{m,n}^1(\alpha)\Lambda_{m,n}^2(\alpha)} = \frac{\vartheta_4(\zeta, q)\vartheta_4(\bar{\zeta}, q)}{P(q)^2} + o(1).$$

4.3. *Computations for $Z_{m,n}^{(01)}(\alpha, \beta)$.* The computations for $Z_{m,n}^{(01)}(\alpha, \beta)$ are slightly different. When performing the product over $w$ for $Z_{m,n}^{(01)}(\alpha, \beta)$, we get

$$Z_{m,n}^{(01)}(\alpha, \beta) = (-1)^{mn} B^{-m} \prod_{j=-\lfloor m/2 \rfloor}^{\lfloor (m-1)/2 \rfloor} \left(1 + Br_1^n\left(\frac{\pi(i\alpha + 2j)}{m}\right)\right)$$

$$\times \left(1 + Br_2^n\left(\frac{\pi(i\alpha + 2j)}{m}\right)\right).$$

In order to obtain an expression similar to what we had above for $Z_{m,n}^{(11)}$ and $Z_{m,n}^{(10)}$, we isolate the term $j = 0$ and factor $r_k^n(\frac{\pi}{m}(i\alpha + 2j))$ when its modulus is greater than 1. This yields

(16)
$$Z_{m,n}^{(01)}(\alpha, \beta)$$
$$= (-1)^{mn}\left(B^{-1} + r_1^n\left(\frac{i\alpha\pi}{m}\right)\right)\left(1 + Br_2^n\left(\frac{i\alpha\pi}{m}\right)\right)$$
$$\times \left(\prod_{j=1}^{\lfloor (m-1)/2 \rfloor} r_1\left(\frac{\pi(i\alpha + 2j)}{m}\right)\right)^n \left(\prod_{j=1}^{\lfloor m/2 \rfloor} r_2\left(\frac{\pi(i\alpha - 2j)}{m}\right)\right)^n$$
$$\times \prod_{j=1}^{\lfloor (m-1)/2 \rfloor} \left(1 + B^{-1}r_1^{-n}\left(\frac{\pi(i\alpha + 2j)}{m}\right)\right)$$

$$\times \left(1 + Br_2^n\left(\frac{\pi(i\alpha + 2j)}{m}\right)\right)$$

(17)
$$\times \prod_{j=1}^{\lfloor m/2 \rfloor} \left(1 + Br_1^n\left(\frac{\pi(i\alpha - 2j)}{m}\right)\right)\left(1 + B^{-1}r_2^{-n}\left(\frac{\pi(i\alpha - 2j)}{m}\right)\right).$$

Let us rewrite the first two terms in brackets, using the fact that $r_1(\frac{i\alpha\pi}{m}) = r_2^{-1}(\frac{i\alpha\pi}{m})$, as

(18)
$$\left(B^{-1} + r_1^n\left(\frac{i\alpha\pi}{m}\right)\right)\left(1 + Br_2^n\left(\frac{i\alpha\pi}{m}\right)\right)$$
$$= \left(B^{-1/2}r_1^{-n/2}\left(\frac{i\alpha\pi}{m}\right) + B^{1/2}r_1^{n/2}\left(\frac{i\alpha\pi}{m}\right)\right)$$
$$\times \left(B^{-1/2}r_2^{-n/2}\left(\frac{i\alpha\pi}{m}\right) + B^{1/2}r_2^{n/2}\left(\frac{i\alpha\pi}{m}\right)\right)$$
$$= 4\cosh\left(\log\left(B^{1/2}r_1^{n/2}\left(\frac{i\alpha\pi}{m}\right)\right)\right)\cosh\left(\log\left(B^{1/2}r_2^{n/2}\left(\frac{i\alpha\pi}{m}\right)\right)\right)$$

and note that

$$\lim_{\substack{n,m\to\infty \\ n/(\sqrt{3}m)\to\rho}} 4\cosh\left(\log\left(B^{1/2}r_1^{n/2}\left(\frac{i\alpha\pi}{m}\right)\right)\right)\cosh\left(\log\left(B^{1/2}r_2^{n/2}\left(\frac{i\alpha\pi}{m}\right)\right)\right)$$
$$= 4\cos\left(\frac{\pi(\alpha\rho - i\beta)}{2}\right)\cos\left(\frac{\pi(\alpha\rho + i\beta)}{2}\right)$$
$$= (2\cos(\bar\zeta))(2\cos(\zeta)).$$

For the other terms of (16), we use the same arguments as for $Z_{m,n}^{(11)}(\alpha, \beta)$ and obtain:

COROLLARY 16. *In the joint limit $m, n \to \infty$, $\frac{n}{\sqrt{3}m} \to \rho$,*

$$\frac{(-1)^{mn}Z_{m,n}^{(01)}(\alpha,\beta)}{\Gamma_{m,n}^1(\alpha)\Gamma_{m,n}^2(\alpha)} = \frac{q^{-1/2}\vartheta_2(\zeta,q)\vartheta_2(\bar\zeta,q)}{P(q)^2} + o(1),$$

*where*

$$\Gamma_{m,n}^1(\alpha) = \left(\prod_{j=1}^{\lfloor(m-1)/2\rfloor} r_1\left(\frac{(i\alpha + 2j)\pi}{m}\right)\right)^n,$$

$$\Gamma_{m,n}^2(\alpha) = \left(\prod_{j=1}^{\lfloor m/2\rfloor} r_2\left(\frac{\pi(i\alpha - 2j)}{m}\right)\right)^n.$$

4.4. *Computations for $Z_{m,n}^{(00)}(\alpha,\beta)$.* Computations are similar to those for $Z_{m,n}^{(01)}(\alpha,\beta)$. The analog of (18) is, in this case,

$$\left(B^{-1} - r_1^n\left(\frac{i\alpha\pi}{m}\right)\right)\left(1 - Br_2^n\left(\frac{i\alpha\pi}{m}\right)\right)$$
$$= 4\sinh\left(\log\left(B^{1/2}r_1^{n/2}\left(\frac{i\alpha\pi}{m}\right)\right)\right)\sinh\left(\log\left(B^{1/2}r_2^{n/2}\left(\frac{i\alpha\pi}{m}\right)\right)\right),$$

which, in the scaling limit, converges to

$$4\sinh\left(\frac{\pi(i\alpha\rho+\beta)}{2}\right)\sinh\left(\frac{\pi(-i\alpha\rho+\beta)}{2}\right)$$
$$= 4\sin\left(\frac{\pi(\alpha\rho-i\beta)}{2}\right)\sin\left(\frac{\pi(\alpha\rho+i\beta)}{2}\right) = (2\sin(\bar\zeta))(2\sin(\zeta)).$$

As in the case of $Z_{m,n}^{(10)}$, the product over the roots of unity introduces a factor $(-1)^{m(n-1)}$. An extra factor $(-1)^{m-1}$ appears since we factor $-r_k(\frac{\pi(i\alpha+2j)}{m})$ when its modulus is greater than 1, something which occurs exactly $\lfloor\frac{m-1}{2}\rfloor + \lfloor\frac{m}{2}\rfloor = m-1$ times. Thus, we have the following result:

COROLLARY 17. *In the joint limit $m, n \to \infty$, $\frac{n}{\sqrt{3}m} \to \rho$,*

$$\frac{(-1)^{mn}Z_{m,n}^{(00)}(\alpha,\beta)}{\Gamma_{m,n}^1(\alpha)\Gamma_{m,n}^2(\alpha)} = -\frac{q^{-1/2}\vartheta_1(\zeta,q)\vartheta_1(\bar\zeta,q)}{P(q)^2} + o(1).$$

4.5. *Investigation of $\Lambda_{m,n}^i(\alpha)$ and $\Gamma_{m,n}^i(\alpha)$, and proof of Proposition 8.* As mentioned in "Ideas governing the next steps of the computations" in Section 4.1, we deal here with that part of each of $Z_{m,n}^{(11)}(\alpha,\beta)$, $Z_{m,n}^{(10)}(\alpha,\beta)$, $Z_{m,n}^{(01)}(\alpha,\beta)$ and $Z_{m,n}^{(00)}(\alpha,\beta)$ which grows exponentially in $mn$. Recall that we expect this growth rate to be driven by the free energy per fundamental domain f of equation (1). Proposition 8 below gives precise statements about this by exhibiting the relation between $\Lambda_{m,n}^1(\alpha)\Lambda_{m,n}^2(\alpha)$, $\Gamma_{m,n}^1(\alpha)\Gamma_{m,n}^2(\alpha)$ and $\Lambda_{m,n}^1(0)\Lambda_{m,n}^2(0)$, and giving the asymptotic behavior of $\Lambda_{m,n}^1(0)\Lambda_{m,n}^2(0)$.

Note that the proof of Proposition 8 is completed by combining Corollaries 14, 15, 16, 17 and Proposition 18.

PROPOSITION 18. *In the joint limit $m, n \to \infty$, $\frac{n}{\sqrt{3}m} \to \rho$,*

1. $A^{n/3} \cdot \frac{\Lambda_{m,n}^1(\alpha)\Lambda_{m,n}^2(\alpha)}{\Lambda_{m,n}^1(0)\Lambda_{m,n}^2(0)} = e^{\pi\alpha^2\rho/2} \cdot (1 + o(1));$
2. $A^{n/3} \cdot \frac{(-1)^n\Gamma_{m,n}^1(\alpha)\Gamma_{m,n}^2(\alpha)}{\Lambda_{m,n}^1(0)\Lambda_{m,n}^2(0)} = q^{1/2} \cdot e^{\pi\alpha^2\rho/2} \cdot (1 + o(1));$
3. $\Lambda_{m,n}^1(0)\Lambda_{m,n}^2(0) = e^{\pi\rho/6} \cdot e^{-mn\mathsf{f}} \cdot (1 + o(1)).$

PROOF. Recall the definition of the functions $r_1(\phi)$ and $r_2(\phi)$ given in (9). Recall, also, that there is an ambiguity in the choice of the square root, which we had clarified on the domain $\overline{D} = \{\phi \in \mathbb{C} | -\pi \leq \mathrm{Re}(\phi) \leq \pi, -c \leq \mathrm{Im}(\phi) \leq c\}$, where $c$ is some positive constant smaller than $\log(4)$.

For the purpose of proving Proposition 18, we need to extend this clarification to a domain $D'$ containing $-\pi \leq \mathrm{Re}(\phi) \leq 3\pi$ and small values of $\mathrm{Im}(\phi)$. On this new domain, we also want to give a meaning to $\log(r_1(\phi))$. Hence, we need to show that we can associate with $r_1(\phi)$ a well-defined argument in $(0, 2\pi)$.

Recall that $r_1(\phi)$ is analytic provided that the term in the square root of (9) does not vanish. This only happens when $(\mathrm{Re}(\phi), \mathrm{Im}(\phi)) = (0[2\pi], -\log(4))$, so the function $r_1$ is analytic on

$$D' = \{\phi \in \mathbb{C} | -\pi < \mathrm{Re}(\phi) < 3\pi, -c' < \mathrm{Im}(\phi) < c'\},$$

where $c'$ is some positive constant smaller than $\log(4)$.

Equating the real parts of the two sides of equation (10) yields

$$\mathrm{Re}(r_1(\phi))\left(1 + \frac{1}{|r_1(\phi)|^2}\right) = e^{-\mathrm{Im}(\phi)} \cos(\mathrm{Re}(\phi)) - 2$$

so that $\mathrm{Re}(r_1(\phi)) < 0$ on $D'$. As a consequence, one can associate, without ambiguity, an argument in $(0, 2\pi)$ [even in $(\frac{\pi}{2}, \frac{3\pi}{2})$] with $r_1(\phi)$. We deduce that the function $\log(r_1)$ is analytic on $D'$.

We need one more property of $r_1$: looking at the term in the square root of (9), we know that $r_1(\phi + 2\pi) = r_2(\phi) = \frac{1}{r_1(\phi)}$. Therefore,

(19) $$\arg(r_1(\phi + 2\pi)) = 2\pi - \arg(r_1(\phi)).$$

*Proof of 1.* Let us first consider the log of the ratio $\frac{\Lambda^1_{m,n}(\alpha)\Lambda^2_{m,n}(\alpha)}{\Lambda^1_{m,n}(0)\Lambda^2_{m,n}(0)}$. Recall that $r_2(\phi) = r_1(\phi + 2\pi)$ so that one can write

$$\Lambda^1_{m,n}(\alpha)\Lambda^2_{m,n}(\alpha) = \left(\prod_{j=0}^{\lfloor(m-1)/2\rfloor} r_1\left(\frac{\pi(i\alpha + 2j + 1)}{m}\right)\right.$$
$$\left.\times \prod_{j=0}^{\lfloor m/2\rfloor - 1} r_1\left(\frac{\pi(i\alpha + 2m - 2j - 1)}{m}\right)\right)^n$$
$$= \left(\prod_{j=0}^{m-1} r_1\left(\frac{\pi(i\alpha + 2j + 1)}{m}\right)\right)^n.$$

It follows that

$$\log \frac{\Lambda^1_{m,n}(\alpha)\Lambda^2_{m,n}(\alpha)}{\Lambda^1_{m,n}(0)\Lambda^2_{m,n}(0)} = n \cdot \sum_{j=0}^{m-1}\left(\log r_1\left(\frac{\pi(i\alpha + 2j + 1)}{m}\right) - \log r_1\left(\frac{\pi(2j + 1)}{m}\right)\right)$$

$$= n \cdot \sum_{j=0}^{m-1} (f_1(\phi_j + \varepsilon_\alpha) - f_1(\phi_j)),$$

where $f_1(\phi) = \log r_1(\phi)$, $\phi_j = \frac{\pi(2j+1)}{m}$ and $\varepsilon_\alpha = \frac{i\pi\alpha}{m}$. We know that $f_1$ is analytic in $D'$, so we can perform a Taylor expansion. This yields

$$\log \frac{\Lambda^1_{m,n}(\alpha)\Lambda^2_{m,n}(\alpha)}{\Lambda^1_{m,n}(0)\Lambda^2_{m,n}(0)}$$

$$= n \cdot \sum_{j=0}^{m-1} \left( \varepsilon_\alpha f_1'(\phi_j) + \frac{\varepsilon_\alpha^2}{2} f_1''(\phi_j) + O\left(\frac{1}{m^3}\right) \right)$$

$$= \frac{i\alpha n}{2} \left( \frac{2\pi}{m} \sum_{j=0}^{m-1} f_1'(\phi_j) \right) - \frac{\alpha^2 \pi n}{4m} \left( \frac{2\pi}{m} \sum_{j=0}^{m-1} f_1''(\phi_j) \right) + O\left(\frac{n}{m^2}\right).$$

For the first term, we have

$$\frac{2\pi}{m} \sum_{j=0}^{m-1} f_1'(\phi_j) = \int_0^{2\pi} f_1'(\phi) \, d\phi - \sum_{j=0}^{m-1} \int_{\phi_j - \pi/m}^{\phi_j + \pi/m} (f_1'(\phi) - f_1'(\phi_j)) \, d\phi$$

$$= f_1(2\pi) - f_1(0) - \sum_{j=0}^{m-1} \int_{\phi_j - \pi/m}^{\phi_j + \pi/m} f_1''(\phi_j)(\phi - \phi_j) \, d\phi + O\left(\frac{1}{m^2}\right)$$

$$= f_1(2\pi) - f_1(0) + O\left(\frac{1}{m^2}\right).$$

For the second term,

$$\frac{2\pi}{m} \sum_{j=0}^{m-1} f_1''(\phi_j) = \int_0^{2\pi} f_1''(\phi) \, d\phi + O\left(\frac{1}{m^2}\right) = f_1'(2\pi) - f_1'(0) + O\left(\frac{1}{m^2}\right).$$

Recombining the different terms gives

$$\log \frac{\Lambda^1_{m,n}(\alpha)\Lambda^2_{m,n}(\alpha)}{\Lambda^1_{m,n}(0)\Lambda^2_{m,n}(0)} = \frac{i\alpha n}{2}(f_1(2\pi) - f_1(0))$$

$$- \frac{\alpha^2 \pi n}{4m}(f_1'(2\pi) - f_1'(0)) + O\left(\frac{n}{m^2}\right),$$

where $f_1(2\pi) = \frac{4i\pi}{3}$, $f_1(0) = \frac{2i\pi}{3}$, $f_1'(2\pi) = -\frac{1}{\sqrt{3}}$, $f_1'(0) = \frac{1}{\sqrt{3}}$.

Hence,

$$\log \frac{\Lambda^1_{m,n}(\alpha)\Lambda^2_{m,n}(\alpha)}{\Lambda^1_{m,n}(0)\Lambda^2_{m,n}(0)} = -\frac{n}{3} \log A + \frac{\alpha^2 \pi \rho}{2} + O\left(\frac{n}{m^2}\right),$$

which proves part 1.

*Proof of 2.* To prove part 2, it suffices to show that $\lim_{\substack{n,m\to\infty \\ n/(\sqrt{3}m)\to\rho}} (-1)^n \Lambda^1_{m,n}(\alpha)$
$\times \Lambda^2_{m,n}(\alpha)/\Gamma^1_{m,n}(\alpha)\Gamma^2_{m,n}(\alpha) = q^{-1/2} = e^{\pi\rho/2}$. We write $\Gamma^1_{m,n}(\alpha)\Gamma^2_{m,n}(\alpha)$ as a unique product:

$$\Gamma^1_{m,n}(\alpha)\Gamma^2_{m,n}(\alpha)$$
$$= \left( \prod_{j=1}^{\lfloor (m-1)/2 \rfloor} r_1\left(\frac{\pi(i\alpha+2j)}{m}\right) \prod_{j=1}^{\lfloor m/2 \rfloor} r_1\left(\frac{\pi(i\alpha+2m-2j)}{m}\right) \right)^n$$

(20)
$$= \left( \prod_{j=1}^{m-1} r_1\left(\frac{\pi(i\alpha+2j)}{m}\right) \right)^n = \left( \prod_{j=1}^{m-1} \sqrt{r_1\left(\frac{\pi(i\alpha+2j)}{m}\right)^2} \right)^n$$

$$= \left( \sqrt{r_1\left(\frac{\pi(i\alpha+2)}{m}\right) r_1\left(\frac{\pi(i\alpha+2(m-1))}{m}\right)} \right.$$
$$\left. \times \prod_{j=1}^{m-2} \sqrt{r_1\left(\frac{\pi(i\alpha+2j)}{m}\right) r_1\left(\frac{\pi(i\alpha+2(j+1))}{m}\right)} \right)^n.$$

The first two terms of the last line consist of the first and last terms of the product (20), while the product over $j$ contains the terms of (20) grouped in pairs, starting from the second one.

From equation (19), we know that $r_1(\phi)r_1(\phi+2\pi) = 1$ and $\arg(r_1(\phi)r_1(\phi+2\pi)) = 2\pi$. In particular, we deduce that

$$\sqrt{r_1\left(\frac{\pi(i\alpha)}{m}\right) r_1\left(\frac{\pi(i\alpha+2m)}{m}\right)} = -1.$$

Therefore,

$$\Gamma^1_{m,n}(\alpha)\Gamma^2_{m,n}(\alpha) = (-1)^n \left( \prod_{j=0}^{m-1} \sqrt{r_1\left(\frac{\pi(i\alpha+2j)}{m}\right) r_1\left(\frac{\pi(i\alpha+2(j+1))}{m}\right)} \right)^n.$$

Let us consider the log of the ratio $\frac{\Lambda^1_{m,n}(\alpha)\Lambda^2_{m,n}(\alpha)}{(-1)^n\Gamma^1_{m,n}(\alpha)\Gamma^2_{m,n}(\alpha)}$. It can be rewritten as

$$\log \frac{\Lambda^1_{m,n}(\alpha)\Lambda^2_{m,n}(\alpha)}{(-1)^n\Gamma^1_{m,n}(\alpha)\Gamma^2_{m,n}(\alpha)}$$
$$= -n \cdot \sum_{j=0}^{m/2-1} \frac{1}{2} \left( f_1\left(\frac{\pi(i\alpha+2j)}{m}\right) - 2f_1\left(\frac{\pi(i\alpha+2j+1)}{m}\right) \right.$$
$$\left. + f_1\left(\frac{\pi(i\alpha+2j+2)}{m}\right) \right).$$

Since the function $f_1$ is analytic in $D'$, a Taylor expansion yields

$$\log \frac{\Lambda^1_{m,n}(\alpha)\Lambda^2_{m,n}(\alpha)}{(-1)^n \Gamma^1_{m,n}(\alpha)\Gamma^2_{m,n}(\alpha)} = -\frac{n\pi}{4m} \cdot \left(\frac{2\pi}{m}\right) \sum_{j=0}^{m-1} f_1''\left(\frac{\pi(2j+1)}{m}\right) + O\left(\frac{n}{m^2}\right)$$

$$= -\frac{n\pi}{4m} \int_0^{2\pi} f_1''(\phi)\, d\phi + O\left(\frac{n}{m^2}\right)$$

$$= -\frac{n\pi}{4m}(f_1'(2\pi) - f_1'(0)) + O\left(\frac{n}{m^2}\right).$$

Hence,

$$\log \frac{\Lambda^1_{m,n}(\alpha)\Lambda^2_{m,n}(\alpha)}{(-1)^n \Gamma^1_{m,n}(\alpha)\Gamma^2_{m,n}(\alpha)} = -\frac{n\pi}{4m}(f_1'(2\pi) - f_1'(0)) + O\left(\frac{n}{m^2}\right)$$

$$= \frac{n\pi}{2\sqrt{3}m} + O\left(\frac{n}{m^2}\right) = \frac{\rho\pi}{2} + O\left(\frac{n}{m^2}\right),$$

which proves part 2.

*Proof of 3.* Recalling the definition of the free energy per fundamental domain f given in equation (1), let us integrate it explicitly over $\psi$:

$$\mathsf{f} = -\frac{1}{4\pi^2} \int_0^{2\pi}\int_0^{2\pi} \log(2(\cos\psi + 1) - e^{i\phi})\, d\phi\, d\psi$$

$$= -\frac{1}{2\pi^2} \int_0^{2\pi}\int_0^{\pi} \log|2(\cos\psi + 1) - e^{i\phi}|\, d\phi\, d\psi \qquad \text{(by symmetry)}$$

$$= -\frac{1}{2\pi^2}\left[\int_0^{\pi}\int_0^{2\pi} \log|e^{i\psi} - r_1(\phi)|\, d\psi\, d\phi + \int_0^{\pi}\int_0^{2\pi} \log|e^{i\psi} - r_2(\phi)|\, d\psi\, d\phi\right]$$

$$= -\frac{1}{\pi}\int_0^{\pi} \log|r_1(\phi)|\, d\phi = \frac{1}{2\pi}\left[\int_0^{\pi} \log r_1(\phi)\, d\phi + \int_0^{\pi} \log r_2(-\phi)\, d\phi\right]$$

$$= -\frac{1}{2\pi}\int_0^{2\pi} f_1(\phi)\, d\phi,$$

where, in the fourth line, we have used the identity

$$\frac{1}{2\pi}\int_0^{2\pi} \log|t + se^{i\psi}|\, d\psi = \begin{cases} \log|t|, & \text{if } |t| \geq |s|, \\ \log|s|, & \text{if } |s| > |t|. \end{cases}$$

Let us now consider the logarithm of $\Lambda^1_{m,n}(0)\Lambda^2_{m,n}(0)$:

$$\log \Lambda^1_{m,n}(0)\Lambda^2_{m,n}(0)$$

$$= n \cdot \sum_{j=0}^{m-1} f_1(\phi_j)$$

$$= \frac{mn}{2\pi}\left[\int_0^{2\pi} f_1(\phi)\,d\phi + \sum_{j=0}^{m-1}\int_{\phi_j-\pi/m}^{\phi_j+\pi/m}(f_1(\phi_j)-f_1(\phi))\,d\phi\right]$$

$$= mn\left[\frac{1}{2\pi}\int_0^{2\pi} f_1(\phi)\,d\phi\right]$$

$$-\frac{\pi n}{12m}(f_1'(2\pi)-f_1'(0))+O\left(\frac{n}{m^3}\right),$$

by the Euler–McLaurin formula. Hence, we deduce that

$$\log \Lambda_{m,n}^1(0)\Lambda_{m,n}^2(0) = -mn\cdot\mathsf{f}+\frac{\pi}{6}\rho+O\left(\frac{n}{m^3}\right). \qquad \square$$

**Acknowledgments.** This work began when the first author was at the Centrum voor Wiskunde en Informatica, Amsterdam. Most of the paper was completed while the second author was at the University of Zürich. The authors would like to thank Bertrand Duplantier, Richard Kenyon and Vincent Pasquier for their helpful comments. We are also grateful to the referees for many pertinent remarks that have helped to improve the quality of this paper.

UPMC UNIVERSITY PARIS 06
UMR 7599
LABORATOIRE DE PROBABILITÉ ET MODÈLES ALÉATOIRES
F-75005, PARIS
FRANCE
E-MAIL: cedric.boutillier@upmc.fr

INSTITUT DE MATHÉMATIQUES
UNIVERSITÉ DE NEUCHÂTEL
RUE EMILE-ARGAND 11
CH-2007 NEUCHÂTEL
SWITZERLAND
E-MAIL: beatrice.detiliere@unine.ch